# Cyclic Difference Sets And Cyclic Hadamard Matrices
N. A. Carella, December 2011.

***Abstract:*** The collection of cyclic Hadamard matrices $\{ H = ( a_{i-j} ) : 0 \leq i, j < n$, and $a_i = -1, 1 \}$ of order $n \geq 1$, is characterized by the orthogonality relation $HH^T = nI$. Only two of such matrices are currently known. It will be shown that this collection consists of precisely two matrices. An application of this result implies that there are exactly seven Barker sequences over the binary set $\{ -1, 1 \}$.

## 1 Introduction

An $n \times n$ matrix $A = ( a_{i,j} )$, $0 \leq i, j < n$, generated by the cyclic shifts of one of its rows, is called a circulant or cyclic matrix and has the form $A = ( a_{i-j} )$, $0 \leq i, j < n$, see [DS] for the analysis of cyclic matrices. The collection of cyclic Hadamard matrices $\{ H = ( a_{i-j} ) : 0 \leq i, j < n$, and $a_i = -1, 1 \}$ of order $n \geq 1$ is characterized by the orthogonality relation $HH^T = nI$. The constraint $HH^T = nI$ calls for the maximal value in the determinant inequality

$$| \det A | \leq \prod_{0 \leq i < n} \| a_{i,j} \|^2 = \prod_{0 \leq i < n} \sum_{0 \leq j < n} | a_{i,j} |^2. \tag{1}$$

In general, Hadamard matrices $H = ( a_{i,j} )$, $0 \leq i, j < n$, can exist for $n = 1, 2$, and $4m$, $m \geq 1$, see [HM] for more details. Currently, there are two cyclic Hadamard matrices $H_1 = < 1 >$ and $H_4 = < 1, 1, 1, -1 >$ known. These are of orders $n = 1$, and $n = 4$ respectively.

This note addresses the following questions stated in [BT, p. 423], and similar sources.

***Conjecture 1.*** There is no cyclic Hadamard matrices of order $v > 4$. Equivalently, there is no cyclic difference sets of parameters $(v, k, \lambda, n) = (4u^2, 2u^2 - u, u^2 - u, u^2)$, $u \geq 1$.

***Conjecture 2.*** There exists no Barker sequences of lengths $v > 13$.

Recent developments in these topics are discussed in [CR], [SC], [MH], [HB], and [CJ]. The information provided within is sufficient to claim that it is plausible to conclude the followings.

***Theorem* 3.**   The collection of cyclic Hadamard matrices has exactly two matrices. These are $H_1 = <1>$ and $H_4 = <1, 1, 1, -1>$.

A Barker sequence $\{s_0, s_1, ..., s_{v-1} : s_i = \pm 1\}$ is a binary sequence of length $v \geq 1$ characterized by the aperiodic autocorrelation function

$$C(\tau) = \sum_{0 \leq i \leq v-1-\tau} s_i s_{i-\tau}, \qquad (2)$$

where $C(\tau) = \pm 1$ for $\tau = 1, 2, ..., v - 1$.

***Theorem* 4.**   The collection of Barker sequences has exactly seven sequences.

The proof of Theorem 3 appears in Subsection 3.1, the proof of Theorem 4 appears in Subsection 3.2.

These problems are studied in great details in [BR], [BT], [BN], [HL], [JB], [TN], and by many other authors. Extensive numerical works, using Number Theoretical methods, has been employed by various authors to eliminate almost every parameters $v < 10^{22}$, see [MH], [SC], et cetera.

## 2 Basic Foundations

This section provides a survey of supporting materials. The interested readers should confer the literature for the detailed analysis of these results.

**Basic Difference Sets**

A difference set $D$ of parameters $(v, k, \lambda, n)$ is a subset $D$ of a group $G$ of orders $k = \#D$ and $v = \#G$ respectively, such that each element $0, 1 \neq d \in G$ has $\lambda \geq 1$ representations as $d = d_1 - d_2, d_i \in D$ in an additive group $G$ or as $d = d_1 d_2^{-1}, d_i \in D$ in a multiplicative group $G$. The integer $n = k - \lambda$ is called the order of the difference set. The difference set is cyclic if the group $G = \mathbb{Z}_v$, the integers modulo $v \geq 1$.

A difference set $D = (v, k, \lambda, n)$ and its complementary difference set $D^c = (v, v - k, v - 2k + \lambda, n)$ occur simultaneously.

The generating function attached to a difference set is defined by the polynomial

$$\theta(X) = \sum_{d \in D} X^d = X^{d_0} + X^{d_1} + X^{d_2} + \cdots + X^{d_{k-1}}. \qquad (3)$$

***Lemma* 5.**   ([BR])   Let $(v, k, \lambda, n)$ be the parameters of a cyclic difference set, and let $w \geq 1$ be a divisor of $v \geq 1$. Then, the generating function satisfies the following congruences.

(i) $\theta_{[w]}(X) = x_0 + x_1 X + x_2 X^2 + \cdots + x_{w-1} X^{w-1} \equiv \theta(X) \bmod X^w - 1$. (4)

(ii) $\theta(X)\theta(X^{-1}) = \sum_{d_i, d_j \in D} X^{d_i - d_j} \equiv n + (\lambda v / w)(1 + X + X^2 + \cdots + X^{w-1}) \bmod X^w - 1$.

The detailed analysis of the derivation of these identities appears in [BR]. The identities for $\theta(X)$ and its various modulo $X^w - 1$ reductions such as $\theta_{[w]}(X)$, and are deployed to generate various relations between the parameters of the cyclic difference set and the coefficients of the generating functions.

**Lemma 6.** ([BR]) Let $(v, k, \lambda, n)$ 0be the parameters of a cyclic difference set, and let $w \geq 1$ be a divisor of $v \geq 1$. Then, the system of Diophantine equations

(1) $x_0 + x_1 + x_2 + \cdots + x_{w-1} = k$, (5)
(2) $x_0^2 + x_1^2 + x_2^2 + \cdots + x_{w-1}^2 = n + \lambda v / w$,
(3) $x_0 x_{0-j} + x_1 x_{1-j} + x_2 x_{2-j} + \cdots + x_{w-1} x_{w-1-j} = \lambda v / w$,

where the index is taken modulo $w$, and $j = 1, 2, \ldots, w - 1$, has a solution $(x_0, x_1, x_2, \ldots, x_{w-1}) \in \mathbb{N}^w$ such that $0 \leq x_i < v / w$ for $i = 0, 1, 2, \ldots, w - 1$.

This is a necessary condition but not sufficient for the existence of a cyclic difference set. Additional, conditions on the variables $x_0, x_1, \ldots, x_{w-1}$ are often used to complete the analysis of existence, for example, the residue class $x_i = \#\{ d \in D : d \equiv i \bmod w \}$ of the variable $x_i$ for $i = 0, 1, 2, \ldots, w - 1$. The c0onverse of this result is not valid. Large scales calculations are given in [GL], and the proof, and related analysis of cyclic difference sets are covered in [BR].

## 3 Cyclic Difference Sets, Hadamard Matrices, and Barker Sequences

The characterization of difference sets by system of Diophantine equations immediately explicates the sporadic nature of the existence of difference sets of certain parameters. The tables of difference sets, see [LJ], clearly confirm this observation on the scarcity of certain difference sets.

### 3.1 Cyclic Difference Set and Cyclic Hadamard Matrices

There are two types of cyclic Hadamard difference sets over the cyclic group $\mathbb{Z}_v = \{ 0, 1, 2, 3, \ldots, v - 1 \}$ discussed in the literature, confer [BR], [BT]. These cyclic difference sets, often called Menon cyclic difference sets, are the followings:

1. The parameters $(v, k, \lambda, n) = (2n - 1, 2n - 1, n - 1, n)$, $n \geq 1$.
2. The parameters $(v, k, \lambda, n) = (4u^2, 2u^2 - u, u^2 - u, u^2)$, $u \geq 0$.

The form of the system of Diophantine equations of the first collection does not preclude the existence of infinite families of cyclic Hadamard difference sets in this collection. But the form of the system of Diophantine equations of the second collection does preclude the existence of infinite families of cyclic Hadamard difference sets because there is a fixed system of equations for all $u \geq 0$, see equations (6), and the parameters are nonlinear functions of $u$.

**The Systems of Equations for $w = 4$.**
The cyclic difference set $(v, k, \lambda, n) = (4u^2, 2u^2 - u, u^2 - u, u^2)$ as a function of $u \geq 0$. The system of Diophantine equations

$$f_0(x_0,...,x_3,u) = x_0 + x_1 + x_2 + x_3 - 2u^2 + u, \qquad f_1(x_0,...,x_3,u) = x_0^2 + x_1^2 + x_2^2 + x_3^2 - u^2 - u^2(u^2 - u), \qquad (6)$$

$$f_2(x_0,...,x_3,u) = x_0 x_3 + x_1 x_0 + x_2 x_1 + x_3 x_2 - u^2(u^2 - u), \qquad f_3(x_0,...,x_3,u) = x_0 x_2 + x_1 x_3 + x_2 x_0 + x_3 x_1 - u^2(u^2 - u),$$

$$f_4(x_0,...,x_3,u) = x_0 x_1 + x_1 x_2 + x_2 x_3 + x_3 x_0 - u^2(u^2 - u),$$

derived from Lemma 6, characterizes cyclic difference set $(v, k, \lambda, n)$ as a function of five variables $x_0 \geq 0$, $x_1 \geq 0$, $x_2 \geq 0$, $x_3 \geq 0$, $u \geq 0$.

The next result will be used to complete the system of equations (6), which has five independent variables and four independent equations. This result is actually an extension of Lemma 6. In Lemma 6 the equations (3) are derived in the ring $\mathbb{Z}[x]/(x^w - 1)$, where $x$ is a transcendental element. But equations (7) are derived in the ring $\mathbb{Z}[\zeta]/(\zeta^w - 1)$, where $\zeta$ is an algebraic element.

**Lemma 7.** Let $(v, k, \lambda, n) = (4u^2, 2u^2 \pm u, u^2 \pm u, u^2)$ be the parameters of a cyclic difference set, and let $w = 4$. Then, the variables $x_0 \geq 0$, $x_1 \geq 0$, $x_2 \geq 0$, $x_3 \geq 0$, and $u \geq 0$ of the system of equations (6) satisfy the followings.

$$f_5(x_0,...,x_3,u) = x_0 x_1 + x_1 x_2 + x_2 x_3 - x_3 x_0 - u^2(u^2 - u), \qquad f_6(x_0,...,x_3,u) = x_0 x_2 + x_1 x_3 - x_2 x_0 - x_3 x_1 - u^2(u^2 - u), \qquad (7)$$

$$f_7(x_0,...,x_3,u) = x_0 x_3 - x_1 x_0 - x_2 x_1 - x_3 x_2 - u^2(u^2 - u), \qquad f_8(x_0,...,x_3,u) = x_0^2 + x_1^2 + x_2^2 + x_3^2 - u^2 - u^2(u^2 - u).$$

Proof: Let $\zeta$ be a primitive 8th root of unity, and consider the generating function congruences

$$\theta(X)\theta(X^{-1}) \equiv n + (\lambda v / w)(1 + X + X^2 + \cdots + X^{w-1}) \bmod X^w - 1, \qquad (8)$$

and

$$\theta_{[w]}(X) = x_0 + x_1 X + x_2 X^2 + \cdots + x_{w-1} X^{w-1} \equiv \theta(X) \bmod X^w - 1, \qquad (9)$$

see Lemma 5. Use (9) to evaluate the congruence (8) at $X = \zeta$, and $X^{-1} = \zeta^{-1}$, and match the coefficients of both congruences to arrive at the followings:

(i) $x_0^2 + x_1^2 + x_2^2 + x_3^2 - u^2 - u^2(u^2 - u)$,     the coefficient of the constant term 1,
(ii) $x_0 x_1 + x_1 x_2 + x_2 x_3 - x_3 x_0 - u^2(u^2 - u)$,     the coefficient of the constant term $\zeta$,
(iii) $x_0 x_2 + x_1 x_3 - x_2 x_0 - x_3 x_1 - u^2(u^2 - u)$,     the coefficient of the constant term $\zeta^2$,
(iv) $x_0 x_3 - x_1 x_0 - x_2 x_1 - x_3 x_2 - u^2(u^2 - u)$,     the coefficient of the constant term $\zeta^3$,

where the reduction table

$$1, \ \zeta, \ \zeta^2, \ \zeta^3, \ \zeta^4 = -1, \ \zeta^{-1} = -\zeta^3, \ \zeta^{-2} = -\zeta^2, \ \zeta^{-3} = -\zeta, \ \text{and} \ \zeta^8 = 1 \qquad (10)$$

was used to simplify the calculations. ∎

**Theorem 3.** There is one cyclic difference set of parameters $(v, k, \lambda, n) = (4u^2, 2u^2 \pm u, u^2 \pm u, u^2)$, $u \geq 0$, in the cyclic group $\mathbb{Z}_v$. In particular, the collection of cyclic Hadamard matrices has exactly two matrices. These are $H_1 = <1>$ and $H_4 = <1, 1, 1, -1>$.

Proof: Apply Lemma 7 to augment the system of Diophantine equations (6) to the complete system of equations

$$f_0 = 0, \quad f_1 = 0, \quad f_2 = 0, \quad f_3 = 0, \quad f_4 = 0, \quad f_5 = 0, \quad f_6 = 0, \quad f_7 = 0 , \tag{11}$$

where variables $0 \leq x_0, x_1, x_2, x_3 \leq v/w = u^2$, and $u \geq 0$ is an integer.

The first term of the Grobner basis of the ideal

$$< f_0, \ f_1, \ f_2, \ f_3, \ f_4, \ f_5, \ f_6, \ f_7 > \ = \ < g_0, \ g_1, \ g_2, \ ..., \ g_{\tilde{m}} > \tag{12}$$

has the explicit expression

$$g_0(x_0,...,x_3,u) = u^4 - u^3 . \tag{13}$$

And its roots are $u = 0$, and 1. Thus, this immediately shows that $u = 1$ is the only nontrivial solution of the system of Diophantine equations (11). ∎

For elementary introduction to Grobner bases, confer [CS], [BW], and for the advanced theory, and practical algorithms, refer to the recent articles in the journals.

Further study of systems of Diophantine equations similar to (11) is an interesting research problem with applications to other related problems in cyclic difference set theory. In the specific case on hand, the algebraic variety

$$V(f_0,...,f_7) = \left\{ (x_0,...,x_4) \in K^5 : f_0(x_0,...,x_4) = 0, ..., f_7(x_0,...,x_4) = 0 \right\}, \tag{14}$$

where $x_4 = u$, and $K$ is an algebraic closure of the rational numbers $\mathbb{Q}$, seems to be of genus $g > 0$. The number of integral points on the intersection $V(f_0) \cap \cdots \cap V(f_7)$ seems to be $\deg f_0 \times \deg f_1 = 8$ points, counting multiplicities. The integral solutions, given as affine points $P = (x_0, x_1, x_2, x_3, x_4)$ on the affine plane $\mathbb{P}^5(\mathbb{Z})$, are

$P_0 = (0, 0, 0, 0, 0)$, multiplicity 3, (15)

$P_1 = (0, 0, 0, 1, 1)$, $P_2 = (0, 0, 1, 0, 1)$,

$P_3 = (0, 1, 0, 0, 1)$, $P_4 = (1, 0, 0, 0, 1)$,

$P_\infty = (\infty, \infty, \infty, \infty, \infty)$, the point at infinity.

The cyclic difference set is parametized by $x_4 = u \geq 0$, each one of these solutions is mapped to the same trivial cyclic difference set of parameters $(v, k, \lambda, n) = (4u^2, 2u^2 - u, u^2 - u, u^2) = (4, 1, 0, 1)$. The corresponding trivial cyclic difference set is represented by $D = \{ 3 \}$, or any of the subset $\{ 0 \}$ or $\{ 1 \}$ or $\{ 2 \}$ or $\{ 3 \} \subset \mathbb{Z}_4 = \{ 0, 1, 2, 3 \}$, the cyclic group of integers modulo 4.

And its nontrivial complement cyclic difference set

$$D^c = (4u^2, 2u^2 + u, u^2 + u, u^2) = (4, 3, 2, 1), \text{ is represented by } D^c = \{ 0, 1, 2 \}.$$

Either one of these cyclic difference sets maps into the circulant Hadamard matrix $H_4 = <1, 1, 1, -1>$, and the Barker sequence $s_v = 1, 1, 1, -1$.

### 3.2 Barker Sequences

Let $D \subset G$ be a difference set in the group $G$. The characteristic function

$$\chi(d) = \begin{cases} 1 & \text{if } d \in G, \text{ and } d \in D, \\ 0 & \text{if } d \in G, \text{ and } d \notin D, \end{cases} \quad (16)$$

of a difference set $D = \{ d_0, d_1, d_2, ..., d_{k-1} \}$ generates a binary sequence $\{ s_0, s_1, s_2, ..., s_{v-1} : s_i = \pm 1 \}$ by setting

$$s_i = (-1)^{1-\chi(d_j)}, \quad d_j \in G, \text{ for } j = 1, 2, ..., v - 1. \quad (17)$$

The periodic autocorrelation has two values:

$$R(\tau) = \sum_{0 \leq i \leq v-1} s_i s_{i-\tau} = \begin{cases} v & \text{if } \tau \equiv 0 \bmod v, \\ v - 4(k - \lambda) & \text{if } \tau \not\equiv 0 \bmod v, \end{cases} \quad (18)$$

for $\tau = 1, 2, ..., v - 1$. The binary sequence attached to this cyclic difference set turns out to be a Barker sequence of aperiodic autocorrelation function $C(\tau) = \sum_{0 \leq i \leq v-1-\tau} s_i s_{i-\tau} = \pm 1$, see equation (2). For example, for the cyclic group of integers modulo 4, $\mathbb{Z}_4 = \{ 0, 1, 2, 3 \}$, there is trivial cyclic difference set $D = \{ 3 \}$, with parameters $(v, k, \lambda, n) = (4, 1, 0, 1)$, and the corresponding Barker sequence is $s_k = 1, 1, 1, -1$.

***Theorem 4.*** The collection of Barker sequences has exactly seven sequences. These are

| | |
|---|---|
| $v = 2$ | $s_v = 1, 1$ |
| $v = 3$ | $s_v = 1, 1, -1$ |
| $v = 4$ | $s_v = 1, 1, 1, -1$ |
| $v = 5$ | $s_v = 1, 1, 1, -1, 1$ |
| $v = 7$ | $s_v = 1, 1, 1, -1, -1, 1-1$ |
| $v = 11$ | $s_v = 1, 1, 1, -1, -1, -1, 1-1, -1, 1, -1$ |
| $v = 13$ | $s_v = 1, 1, 1, 1, 1, -1, -1, 1, 1, -1, 1, -1, 1.$ |

Proof: The case $v =$ odd is due to Turyn and Storer in [TN], and the case $v =$ even follows from Theorem 3.  ∎

Except for a few cases, the nonexistence of Barker arrays in dimension $n \geq 2$ has also been proven, see [JM].

The complete collection of nontrivial cyclic difference sets associated with Baker sequences are listed below, the map (17) maps the cyclic difference set $D = (v, k, \lambda, n)$ to the binary sequence $\{ s_0, s_1, s_2, ..., s_{v-1} : s_i = \pm 1 \}$.

$v = 3$, the parameters are $(3, 2, 1, 1)$ and the cyclic difference set is $D = \{ 0, 1 \}$.

$v = 4$, the parameters are $(4, 3, 2, 1)$ and the cyclic difference set is $D = \{ 0, 1, 2 \}$.

$v = 5$, the parameters are $(5, 4, 3, 1)$ and the cyclic difference set is $D = \{ 0, 1, 2, 4 \}$.

$v = 7$, the parameters are $(7, 4, 2, 2)$ and the cyclic difference set is $D = \{ 0, 1, 2, 5 \}$.

$v = 11$, the parameters are $(11, 5, 2, 3)$ and the cyclic difference set is $D = \{ 0, 1, 2, 6, 9 \}$.

$v = 13$, the parameters are $(13, 9, 6, 3)$ and the cyclic difference set is $D = \{ 0, 1, 2, 3, 4, 7, 8, 10, 12 \}$.

# 4 References

[BN] P. Borwein, Computational excursions in analysis and number theory. CMS Books in Mathematics/Ouvrages de Mathématiques de la SMC, 10. Springer-Verlag, New York, 2002.

[BR] Baumert, Leonard D. Cyclic difference sets. Lecture Notes in Mathematics, Vol. 182 Springer-Verlag, Berlin-New York 1971.

[BT] Beth, Thomas; Jungnickel, Dieter; Lenz, Hanfried Design theory. Vol. I. Second edition. Encyclopedia of Mathematics and its Applications, 69. Cambridge University Press, Cambridge, 1999. ISBN: 0-521-44432-2.

[BW] Becker, Thomas; Weispfenning, Volker Gröbner bases. A computational approach to commutative algebra. In cooperation with Heinz Kredel. Graduate Texts in Mathematics, 141. Springer-Verlag, New York, 1993. ISBN: 0-387-97971-9.

[CR] Craigen, R.; Kharaghani, H. On the nonexistence of Hermitian circulant complex Hadamard matrices. Australas. J. Combin. 7 (1993), 225–227.

[CJ] Robert Craigen, Jonathan Jedwab, Comment on "The Hadamard circulant conjecture", 10 February 2011,

[CS] Cox, David; Little, John; O'Shea, Donal Ideals, varieties, and algorithms. An introduction to computational algebraic geometry and commutative algebra. Third edition. Springer, New York, 2007. ISBN: 978-0-387-35650-1; 0-387-35650-9.

[DS] Davis, P. J. Circulant Matrices, 2nd ed. New York: Chelsea, 1994.

[GL] Gaal, Peter; Golomb, Solomon W. Exhaustive determination of $(1023,511,255)$(1023,511,255) -cyclic difference sets. Math. Comp. 70 (2001), no. 233, 357–366.

[HB] B. Hurley, P. Hurley, T. Hurley, The Hadamard circulant conjecture, January 24, 2011, Bull. London Math. Soc. 11 Pages.

[HL] T. Helleseth, P. Vijay Kumar. Sequences with low correlation, Handbook of Coding Theory, North-Holland, 1998, pp. 1765-1855.

[HM], Horadam, K. J. Hadamard matrices and their applications. Princeton University Press, Princeton, NJ, 2007. ISBN: 978-0-691-11921-2; 0-691-11921-X.

[JB] J. Jedwab, A survey of the merit factor problem for binary sequences, to appear, in T. Helleseth et al, eds., Lecture Notes in Computer Science, Sequences and Their Applications - Proceedings of SETA '04, Springer-Verlag, Preprint 2005.

[JM] Jedwab, Jonathan; Parker, Matthew G. There are no Barker arrays having more than two dimensions. Des. Codes Cryptogr. 43 (2007), no. 2-3, 79–84.

[LJ] La Jolla Cyclic Difference Set Repository, www.ccrwest.org/diffsets.html

[MH] Mossinghoff, Michael J. Wieferich pairs and Barker sequences. Des. Codes Cryptogr. 53 (2009), no. 3, 149–163.

[MT] S. Merten, Exhaustive Search for Low Autocorrelation Binary Sequences, arXiv:cond-mat/9605050v2, Preprint 1996.

[MW] mathworld.com

[SC] Ka Hin Leung, Bernhard Schmidt, New restrictions on possible orders of circulant Hadamard matrices, Preprint, 2011, Designs, Codes and Cryptography.

[TN] R. Turyn, J. Storer. On binary sequences. Proc. Amer. Math. Soc. 12 1961, pp. 394-399.

[WK] en.wikipedia.org.